\documentclass[11pt]{article}
\usepackage{amssymb}
\usepackage{latexsym,amsfonts,amsmath}
\usepackage{color}

\normalbaselines

\newtheorem{Theorem}{Theorem}[section]

\newtheorem{Definition}[Theorem]{Definition}

\newtheorem{Proposition}[Theorem]{Proposition}

\newtheorem{Lemma}[Theorem]{Lemma}

\newtheorem{Corollary}[Theorem]{Corollary}

\newtheorem{Remark}[Theorem]{Remark}

\newcommand{\RR}{{{\rm I} \kern -.15em {\rm R} }}

\newcommand{\C}{{{\rm l} \kern -.42em {\rm C} }}

\newcommand{\nat}{{{\rm I} \kern -.15em {\rm N} }}

 \newcommand{\ds}{\displaystyle}

\newcommand{\be}{\begin{equation}}
\newcommand{\ee}{\end{equation}}
\newcommand{\beq}{\begin{eqnarray}}
\newcommand{\eeq}{\end{eqnarray}}
\newcommand{\beqs}{\begin{eqnarray*}}
\newcommand{\eeqs}{\end{eqnarray*}}
\newcommand{\bt}{\begin{Theorem}}
\newcommand{\et}{\end{Theorem}}
\newcommand{\br}{\begin{Remark}}
\newcommand{\er}{\end{Remark}}
\newcommand{\bc}{\begin{Corollary}}
\newcommand{\ec}{\end{Corollary}}
\newcommand{\bl}{\begin{Lemma}}
\newcommand{\el}{\end{Lemma}}
\newcommand{\bd}{\begin{definition}}
\newcommand{\ed}{\end{definition}}

\title{Stability of solutions to nonlinear wave
 equations\\
with
switching time--delay}
\author{{\sc Genni Fragnelli}
\\Dipartimento di Matematica\\Universit\`a di Bari\\
Via E. Orabona 4, 70125 Bari, Italy
\\\\
{\sc Cristina Pignotti}
\\Dipartimento di Ingegneria e Scienze dell'Informazione e Matematica\\
 Universit\`{a} di L'Aquila\\
Via Vetoio, Loc. Coppito, 67010 L'Aquila, Italy}

\date{}

\begin{document}

\textwidth=160 mm

\textheight=225mm

\parindent=8mm

\frenchspacing

\maketitle

\begin{abstract}
In this paper we study well--posedness and asymptotic stability for a class of nonlinear second-order evolution equations with intermittent delay damping. More precisely,
a delay feedback and an undelayed one act alternately in time. We show that, under suitable conditions on the feedback operators, asymptotic stability results are available. Concrete examples included in our setting are illustrated.
We give also stability results for an abstract model with alternate  positive--negative damping, without delay.
\end{abstract}

\vspace{5 mm}

\def\qed{\hbox{\hskip 6pt\vrule width6pt
height7pt
depth1pt  \hskip1pt}\bigskip}



\section{Introduction}
\label{pbform}\hspace{5mm}

\setcounter{equation}{0}

Let $H$ be a real Hilbert space with scalar product $\langle\cdot, \cdot\rangle_H$  and norm $\Vert \cdot\Vert ,$ and let $A:{\mathcal D}(A)\rightarrow H$
be a linear self--adjoint  coercive operator on $H$ with dense domain. Let $V:={\mathcal D}(A^{\frac 1 2})$, the domain of
$A^{\frac 1 2}$ with norm $\Vert v\Vert_V=\Vert A^{\frac 12}v\Vert_H,$ be such
that
$$V\hookrightarrow H\equiv H^\prime  \hookrightarrow V^\prime,$$
with dense embeddings. Then, there exists $\lambda_1>0$ such that
\begin{equation}\label{autoval}
\lambda_1\Vert u\Vert_H^2\le \Vert u\Vert_V^2, \quad \forall \ u\in  V.
\end{equation}

Moreover,
let $U_i,\ i=1,2,$ be
     real Hilbert spaces
with norm and inner product denoted respectively by $\Vert \cdot\Vert_{U_i}$
and $\langle\cdot ,\cdot\rangle_{U_i}$
and let
$B_i(t)\in {\mathcal L} (U_i,  {H}),$ $i=1,2,$ be time--dependent operators satisfying

$$B_1^*(t)B_2^*(t)=0,\quad \forall\; t>0.$$

Let us consider the problem
 \begin{eqnarray}
& &u_{tt}(t) +A u (t)+B_1(t)B_1^*(t)u_t(t) +B_2(t)B_2^*(t)u_t(t-\tau) =f(u),\quad t>0,\label{1.1}\\
& &u(0)=u_0 \in V\quad \mbox{\rm and}\quad u_t(0)=u_1\in H,\quad\label{1.2}
\end{eqnarray}
where the constant  $\tau >0$ is the time delay.

The example we have in mind is, for $p\ge 0,$
 \begin{eqnarray}
& &u_{tt}(x,t) -\Delta u (x,t)+b_1(t) u_t(x,t)+b_2(t)  u_t(x,t-\tau)=-\vert u\vert^p u\ \mbox{\rm in}\ \Omega\times
(0,+\infty),\label{Wo.1}\\
& &u (x,t) =0\quad \mbox{\rm on}\quad\partial\Omega\times
(0,+\infty),\label{Wo.2}\\
& &u(x,0)=u_0(x)\quad \mbox{\rm and}\quad u_t(x,0)=u_1(x)\quad \hbox{\rm
in}\quad\Omega,\label{Wo.3}
\end{eqnarray}
with
initial
data $(u_0, u_1)\in H^1_0(\Omega)\times L^2(\Omega)$, where $\Omega$ is a bounded and smooth domain of $\RR^N$
and
$b_1,b_2$ in $L^\infty (0,+\infty)$ are such that
$$b_1(t)b_2(t)=0,\quad\forall\ t>0.$$

In this case $H= U_i= L^2(\Omega)$, $B_i^*(\varphi)= \sqrt{b_i}\varphi$ for all $\varphi \in H$ and $i=1,2$, $V= H^1_0(\Omega)$  and $\lambda_1$ in (\ref{autoval}) is the first eigenvalue of $-\Delta$ on $H^1_0(\Omega)$, being (\ref{autoval}) the usual  Poincar\'e's inequality.

Time delay is often present
in  applications and practical problems and it is by now well--known that even an arbitrarily small delay in the feedback may destabilize a system wich is uniformly exponentially stable in absence of delay.
See, e.g.,  \cite{Datko, DLP, NPSicon06, XYL} where examples of the destabilization effect due to time delay are given.

The idea is then to use a stabilizing feedback in order to contrast the instability due to the presence of a delay term.
In \cite{NPSicon06, XYL} a standard damping and a delayed one act simultaneously and the stability of the systems is guaranteed if the coefficient of the undelayed damping is bigger than the one of the delay feedback.
Then, in \cite{ADE2012, NP2014}, the authors consider the case of delayed--undelayed feedback acting in alternate time intervals and give sufficient conditions on the feedback operators in order to have stability.

A similar problem has been
considered
        in \cite{ANP2012} for the one dimensional wave equation but with a  completely           different strategy.
Indeed,
  in \cite{ANP2012} stability results are obtained, only for particular values of the time delays related to the length of the domain (cfr. \cite{Gugat}), by using the D'Alembert formula.

Here we extend the results of \cite{ADE2012, NP2014} to nonlinear models but also we
significantly improve some results there proved. Indeed, we are able to remove an assumption (see (\ref{new}) below)  on the feedback bounds obtaining more general results also in the linear case.

 By using suitable observability inequalities  for the  model with only the undelayed feedback  and through the definition of a suitable energy (see (\ref{energy2A})), we obtain
sufficient conditions ensuring asymptotic stability.

Some concrete examples falling in our abstract setting are also illustrated.

With the same approach we also consider nonlinear second order evolution equations without delay but with positive--negative dampings acting alternately.
This kind of problem was first considered in \cite{HMV} in the linear case and then extended to nonlinear models in \cite{Genni2}.
In both papers only the case of distributed damping was considered. Here, we meaningfully generalize these results by considering the case of local damping. More precisely, in concrete examples, the positive (stabilizing) damping and the negative (destabilizing) one
may be localized in whatever subsets of the domain. The only geometric requirement is, of course, that the positive damping has to be localized in a region satisfying a control geometric property (see \cite{BLR}).

The interest for models with intermittent delayed--undelayed damping or positive--negative dampings is motivated by various applications.  For example, the presence of  positive--negative damping can be found in aerodynamics:
nose wheel shimmy of an airplane is the consequence of a negative
damping, which is controlled by a suitable hydraulic shimmy damper
which induces a positive damping (\cite{S}). Another example of
sign--changing damping comes from Quantum Field Theory and Landau
instability (see \cite{KN}) and from mesodynamics with the laser
driven pendulum (see\cite{ds}). Actually, negative damping may appear in
every--day--life, for example Gunn diodes, used as source of
microwave power, and suspension bridges (\cite{LaMK}, \cite{bila},
\cite{mm1}, \cite{mm2}), which may experience negative damping in a
catastrophic way, like Takoma Bridge.
Observe also that the recent results given in \cite{GGH} show that  dampings with pulsating coefficients are more effective, with respect to the ones with constant coefficients, in order to stabilize second order evolution equations.
This is a further motivation for our study.

The paper is organized as follows. In section \ref{well} we introduce our abstract setting and give a well--posedness result. In section \ref{st} we prove the asymptotic stability results. We consider first distributed dampings, then the localized case and finally, for the linear model, we give the results under more explicit conditions.
Finally, in section \ref{PN} we consider the model without delay and positive--negative dampings.

\section{The abstract setting\label{well}}

\hspace{5mm}

\setcounter{equation}{0}

In order to deal with the well-posedness of
\eqref{1.1} -- \eqref{1.2}, first
we consider  the abstract problem
\[
u_{tt}(t)+ Au(t) + B(t)u_t(t)=f(u)
\]
and its associated Cauchy problem
\begin{eqnarray}
& & u_{tt}(t)+ Au(t) + B(t)u_t(t)=f(u) ,\quad t>0, \label{wave}\\
& & u(0)=u_0\in V\quad \mbox{\rm and}\quad  u_t(0)=u_1\in H. \label{wave1}
\end{eqnarray}
Here $H$ and $V$ are as before and $B= B_1B_1^*: V \rightarrow V^\prime.$ We recall  the next definition:
\begin{Definition}
A function $u$ is a weak solution of \eqref{wave} -- \eqref{wave1} if
for any $T>0$ we have  $$u\in L^2(0,T;V)\cap H^1(0,T;H)\cap
H^2(0,T;V')
$$
with $B(t)u_t(t)\in H$ for any $t$, $\langle Bu_t,u_t\rangle_H \in
L^2(0,T)$ and
\[
Au\in L^2(0,T;V'),\quad  Bu_t\in L^2(0,T;V'),\quad f(u)\in
L^2(0,T;H).
\]
Moreover, $u$ is such that
$u(0)=u_0,\ u_t(0)=u_1$ and
\[
u_{tt}(t)+ Au(t)+B(t)u_t(t)=f(u) \quad \mbox{ in }L^2(0,T;V').
\]
\end{Definition}
Now, rewrite \eqref{wave}  -- \eqref{wave1} as
\begin{eqnarray}
& & U_t + LU + C(U) =0,\label{wave2}\\
& & U(0)= U^0, \label{wave3}
\end{eqnarray}
where $U=(u, u_t)$,  $L = \begin{pmatrix}  0 &-I\\
A& B\end{pmatrix} $, $C( U):= \begin{pmatrix} 0\\
-f(u)\end{pmatrix},$
and $U^0=(u_0, u_1)$.

On the nonlinear term $f$ we assume
\begin{equation}\label{l.lip.}
\begin{aligned}
f \;& \text{is  locally Lipschitz continuous, i.e. }\\
& \forall \ K>0\ \ \exists \ L(K)\ \ \mbox{\rm such that}\ \ \|f(u) - f(v)\|_H \le L(K) \|u- v\|_V,
\end{aligned}
\end{equation}
provided $\|u\|_V, \|v\|_V \le K$;
\begin{equation}\label{NL1}
sf(s)\le 0,\quad \forall\  s\in\RR,
\end{equation}
which implies
$$F(s):=\int_0^s f(r) dr\le 0,\quad \forall\ s\in \RR\,$$
or
\begin{equation}\label{NL2}
sf(s)-F(s)\le 0, \quad\forall\ s\in\RR\,.
\end{equation}

As prototype, we can consider the function $f(u)=-|u|^pu,$ $p\ge0$. Clearly $f$ is locally Lipschitz continuous. Moreover, we remark that the sign assumptions on $f$ is quite reasonable and hard to to relax. Indeed, Levin, Park and Serrin in \cite{l} and \cite{lps} proved that the solutions of $u_{tt} -\Delta u+ a(x,t)u_t= |u|^p u$ in $\Omega$ with $p > 0$ and $a(x,t) \ge 0$ can blow up in finite time.

Observe that, setting $\mathcal H:= V \times H,$ \eqref{l.lip.} implies that $C: \mathcal H \rightarrow \mathcal H$ is locally Lipschitz continuous. Hence, defining
\[
D(L):=\{ (u,v) \in V \times V: Au + Bv \in H \},
\]
 we can apply \cite[Theorem 7.2]{lasiecka}, obtaining the following existence result:
\begin{Theorem}\label{WP}
Suppose that $L$ is a maximal monotone mapping, $L0=0$ and $U^0\in D(L)$.
Then there exists $T_M$ such that problem \eqref{wave2} -- \eqref{wave3} has a
unique strong solution $U$ on the interval $[0, T_M)$, i.e. $U \in W^{1,\infty}(0, T_M; \mathcal H)$.
Furthermore, if we assume only $U^0 \in \mathcal H$ we obtain a unique weak solution $U \in C([0, T_M); \mathcal H)$.
\\
In both cases we have
\[
\lim_{t \rightarrow T_M} \|u(t)\|_V = \infty,
\]
provided $T_M < \infty$.
\end{Theorem}
Observe that if $A$ is a self-adjoint, positive and coercive operator with dense domain in $H$  and if $B \in \mathcal L(V,V^\prime)$ is such that $\langle Bv, v \rangle_H \ge 0$ for all $v \in V$, then $L$ is a maximal monotone operator with dense domain
in $\mathcal H$ (see \cite{aloui}, \cite{aloui1}). Hence, as a consequence of Theorem \ref{WP}, we have:
\begin{Corollary}\label{WP1}
Assume that $A$ is a self-adjoint, positive and coercive operator with dense domain in $H$  and $B \in \mathcal L (V,V^\prime)$ is such that $\langle Bv, v \rangle_H \ge 0$ for all $v \in V$. If $(u_0, u_1) \in D(L)$ then there exists $T_M$ such that problem \eqref{wave} -- \eqref{wave1} has a
unique strong solution $u$ on the interval $[0, T_M)$, i.e. $u \in W^{1,\infty}(0, T_M; V)$.
Furthermore, if we assume only $(u_0, u_1) \in \mathcal H$ we obtain a unique weak solution $(u, u_t) \in C([0, T_M); \mathcal H)$.
\\
In both cases we have
\[
\lim_{t \rightarrow T_M} \|u(t)\|_V = \infty,
\]
provided $T_M < \infty$.
\end{Corollary}

Now, for any solution of problem \eqref{wave} -- \eqref{wave1}, we consider the energy associated to such a solution:
 \begin{equation}\label{energy standard}
 \displaystyle{
E_S(t)=E_S(u;t):= \frac{1}{2}\Big(\|u(t)\|_V^2+\|u_t(t)\|_H^2\Big )
 -{\mathcal F}(u),}
 \end{equation}
 where $\mathcal F$ is a real-valued functional such that $\mathcal F(0)=0$ and $\mathcal F'(u) (v) = \langle f(u), v \rangle_{V',V}$ for all $u, v \in V$. Of course, in problem (\ref{Wo.1}) -- (\ref{Wo.3})
$${\mathcal F}(u)=\int_{\Omega} F(u) dx,$$
where $F(s)=\int_0^s f(t) dt,$ i.e. $F(s)=-\ds\frac {\vert s\vert^{p+2}}{p+2}$ for the model case.
The following existence result holds
\begin{Theorem}\label{WP2}
Assume that $A$ is a self-adjoint, positive and coercive operator with dense domain in $H$, $B \in \mathcal L (V,V^\prime)$ is such that $\langle Bv, v \rangle_H \ge 0$ for all $v \in V$ and  $\mathcal F \le 0$. Moreover, assume that there exists a positive constant $C$ such that $E_S(T) \le C E_S(0)$ for all $T \in (0, T_M)$. If $(u_0, u_1) \in D(L)$ then problem \eqref{wave} -- \eqref{wave1} has a
unique strong solution $u$ on the interval $[0, \infty)$.
Furthermore, if $(u_0, u_1) \in \mathcal H$ we obtain a unique weak solution $(u, u_t) \in C([0, \infty); \mathcal H)$.
\end{Theorem}
\noindent {\bf Proof.} Thanks to Corollary \ref{WP1}, we know that there exists a unique solution in $[0, T_M)$. Assume, by contradiction, that $T_M < \infty$. Then
\\
\begin{equation}\label{BU}
\lim_{t \rightarrow T_M} \|u(t)\|_V = \infty.
\end{equation}
By definition of $E_S(t)$ and since $\mathcal F \le 0$, it follows that
\[
\|u(t)\|_V^2  \le 2 E_S(T) \le 2 C E_S(0).
\]
Hence \eqref{BU} cannot happen.
\qed

Clearly, if  \eqref{NL1} is satisfied, then $\mathcal F \le 0$. Moreover, observe that in the linear case, i.e. $f \equiv 0$, the existence and uniqueness of a solution in $[0, \infty)$ is guaranteed, for example,  by \cite[Theorem 7.1]{lasiecka}.

\bigskip

Now, we assume that for all $n\in\nat$, there exists $t_n>0$, with $t_n<t_{n+1}$,
such that
\begin{eqnarray*}
B_2(t)=0\  \forall\ t\in I_{2n}=[t_{2n},t_{2n+1}),\\
B_1(t)=0
\ \forall \  t\in I_{2n+1}=[t_{2n+1},t_{2n+2}),
\end{eqnarray*}
with $B_1\in C^1([t_{2n},t_{2n+1}]; {\mathcal L} (U_1,  {H}))$
and $B_2\in C^1([t_{2n+1},t_{2n+2}]; {\mathcal L} (U_2,  {H}))$.
We further assume
 \begin{equation}\label{length}
 \tau\leq T_{2n},\quad \forall \  n\in\nat,
  \end{equation}
 where $T_n$ denotes the length of the interval $I_n,$ that is
\begin{equation}\label{Tn}
T_n=t_{n+1}-t_n,\quad n\in \nat.
\end{equation}

Let $W$ be an Hilbert space such that $H$ is continuously embedded into $W,$ i.e.
\begin{equation}\label{embedding}
  \|u\|_W^2\le C \|u\|_H^2,\quad\forall \;u\in H \ \mbox{\rm with}\ \  C>0\ \
\mbox{\rm  independent of}\  u.
\end{equation}

We assume that, for all $n\in\nat$, there exist three positive constants $m_{2n}$,
 $M_{2n}$ and $M_{2n+1}$, with $m_{2n}\leq M_{2n}$, such that for all $u\in H$ we have

\vskip+5pt

i) $m_{2n}\|u\|_W^2\le  \|B_1^*(t)u\|_{U_1}^2\le M_{2n} \|u\|_W^2$ for $t\in I_{2n}=[t_{2n},t_{2n+1}),$  $\ \forall\ n\in\nat;$

\vskip+2pt

ii)$\|B_2^*(t)u\|_{U_2}^2 \le M_{2n+1}\|u\|_W^2 $ for $t\in I_{2n+1}=[t_{2n+1},t_{2n+2}),$ $\ \forall\ n\in\nat.$

\vskip+5pt

Let us introduce the energy functional
\begin{equation}\label{energy2A}
 \displaystyle{
E(t)=E(u;t):= \frac{1}{2}\Big(\|u(t)\|_V^2+\|u_t(t)\|_H^2\Big ) +
\frac{1} 2
 \int_{t-\tau}^{t}\Vert B^*_2(s+\tau)u_t(s)\Vert_
{U_2}
^2  ds -{\mathcal F}(u).
}
\end{equation}

Note that
(\ref{energy2A}) is the usual energy $E_S(\cdot),$  for wave-type equation in presence of the nonlinearity $f,$
plus an integral term (see \cite{ADE2012}, cfr. also \cite{NPSicon06}) due to the presence of the time delay.

Now, we give the following definition:
\begin{Definition}
A solution of problem $(\ref{1.1})-(\ref{1.2})$  is a function $u$ such that for any $T>0$
$$u\in L^2(0,T;V)\cap H^1(0,T;H)\cap H^2(0,T; V^\prime )$$
with $\Vert B_1^*u_t\Vert_{U_1}\in L^2(0,T),$ $\Vert B_2^*u_t(\cdot -\tau)\Vert_{U_2}\in L^2(0,T)$  and
$$Au \in L^2(0,T;V^\prime),\  B_1^*u_t\in L^2(0,T;V^\prime ),\ B_2^*u_t(\cdot-\tau)\in L^2(0,T;V^\prime ), \ f(u)\in L^2(0,T; H).$$
Moreover, $u$ is such that $u(0)=u_0,$ $u_t(0)=u_1$ and
$$u_{tt}(t) +A u (t)+B_1(t)B_1^*(t)u_t(t) +B_2(t)B_2^*(t)u_t(t-\tau) =f(u) \ \mbox{\rm in}\ L^2(0,T; V^\prime ).$$
\end{Definition}

Observe that if  $H=U_i= L^2(\Omega)$, $i=1,2$, and $V= H^1_0(\Omega)$, then the condition $f(u) = -|u|^p u \in L^2(0,T;H)$ is clearly satisfied when $p \ge 0$ if $N=1,2$ or $0< p \le \frac{2}{N-2}$ if $N \ge 3$.
\vspace{0.5cm}

\begin{Remark}{\rm
Our assumptions do not ensure
that the energy $E(\cdot )$  is  decreasing on the time intervals $I_{2n}$ where only the standard frictional damping acts, i.e. $B_2\equiv 0\,,$ as of course it happens for the standard energy $E_S(\cdot).$
In order to have a decay estimate for $E(\cdot )$ in the intervals $I_{2n},$
we should assume, as in \cite{NP2014},
\begin{equation}\label{new}
\inf_{n\in\nat}
 \frac {m_{2n}}{M_{2n+1}} >0,
\end{equation}
and  define $E(\cdot)$ as
$$
 \displaystyle{
E(t)=E(u;t):= \frac{1}{2}\Big(\|u(t)\|_V^2+\|u_t(t)\|_H^2\Big ) +
\frac{\xi} 2
 \int_{t-\tau}^{t}\Vert B^*_2(s+\tau)u_t(s)\Vert_
{U_2}
^2  ds -{\mathcal F}(u).
}
$$
where $\xi$ is a positive number satisfying
$$
\xi <  \inf_{n\in\nat}
\frac {m_{2n}}{M_{2n+1}}.
$$
However, here we do not need $E$ decreasing in the time intervals without delay $I_{2n},$ since in these time intervals we will work with the standard energy $E_S
(\cdot ).$ Consequently,  we do not assume (\ref{new}) to obtain our stability results.
}\end{Remark}

\begin{Proposition}\label{derivE2abstrait}
Assume $\mbox{\rm i),\ ii)},$ $(\ref{NL1})$ and $(\ref{length}).$
For any regular solution of problem $(\ref{1.1})-(\ref{1.2}),$ the energy $E(t)$ satisfies
\begin{equation}\label{stimaderD2abstrait}
E^{\prime}(t)\le
{M _{2n+1}}
\|u_t\|_W^2,
\end{equation}
 for $t\in I_{2n+1},$ $n\in\nat.$
\end{Proposition}

\noindent{\bf Proof:} Differentiating the energy functional, we have
$$E^\prime (t)= \langle u_t,u\rangle_V+\langle u_{tt}, u_t\rangle_H +
\frac {1} 2\Vert B_2^*(t+\tau )u_t(t)\Vert_{U_2}^2 -
\frac {1} 2\Vert B_2^*(t)u_t(t-\tau )\Vert_{U_2}^2-\langle f(u), u_t\rangle_H.$$
Then, from equation (\ref{1.1}),
$$\begin{array}{l}
\displaystyle{ E^\prime (t)=\langle u_t, u_{tt}+Au-f(u)\rangle_{V,V^\prime} +\frac {1} 2\Vert B_2^*(t+\tau )u_t(t)\Vert_{U_2}^2 -
\frac {1} 2\Vert B_2^*(t)u_t(t-\tau )\Vert_{U_2}^2}\\\medskip
\hspace{1 cm}\displaystyle{
=-\langle u_t, B_1(t)B_1^*(t)u_t(t) +B_2(t)B_2^*(t)u_t(t-\tau) \rangle_{V,V'}}\\\medskip
\hspace{1 cm}
\displaystyle{
+\frac{1}{2}\Vert B^*_2(t+\tau)u_t(t)\Vert_{U_2}^2
-\frac{1}{2}\Vert B^*_2(t)u_t(t-\tau)\Vert_{U_2}^2}.
\end{array}
$$
Therefore we obtain
$$
\begin{array}{l}
\displaystyle{
E^\prime(t)=-\|B_1^*(t)u_t(t)\|^2_{U_1} -\langle B_2^*(t)u_t,B_2^*(t)u_t(t-\tau)\rangle_{U_2}}\\\medskip
\hspace{1 cm}\displaystyle{
+\frac{1}{2}\Vert B^*_2(t+\tau)u_t(t)\Vert_{U_2}^2
-\frac{1}{2}\Vert B^*_2(t)u_t(t-\tau)\Vert_{U_2}^2.}
\end{array}
$$
For $t\in I_{2n+1},$ it is $B_1(t)=0$ and so  the previous identity gives
\begin{eqnarray*}
E^\prime (t)=-\langle B_2^*(t)u_t,B_2^*(t)u_t(t-\tau)\rangle_{U_2}+
\frac{1}{2}\Vert B^*_2(t+\tau)u_t(t)\Vert_{U_2}^2
-\frac{1}{2}\Vert B^*_2(t)u_t(t-\tau)\Vert_{U_2}^2.
\end{eqnarray*}
By using Young's inequality we have
\begin{eqnarray*}
E^\prime(t)\le\frac{1}{2}\Vert B^*_2(t)u_t(t)\Vert_{U_2}^2
+\frac{1}{2}\Vert B^*_2(t+\tau)u_t(t)\Vert_{U_2}^2.
\end{eqnarray*}
This proves (\ref{stimaderD2abstrait}) using assumption ii) because $t+\tau$ belongs either to $I_{2n+1},$
or to $I_{2n+2}$ and in the last case $B^*_2(t+\tau)=0$.
$\qed$

Proceeding analogously to \cite{ADE2012} and using Theorem \ref{WP2} we can prove the following existence result.
\begin{Theorem}\label{WP3}
Under the assumptions of Theorem \ref{WP2}, if $(u_0, u_1) \in V \times H$, for any $T >0$ we obtain a unique weak solution
$$u\in C([0,T];V)\cap C^1([0,T];H).$$
\end{Theorem}
\noindent {\bf Proof.} We can combine analogous lemma in  \cite{GenniSicon} with the well--posedness result in \cite{ADE2012}.
We can argue on the interval $[0,t_2)$ which is the union of the first time interval $[0,t_1),$ where the delay term is no present, and the  second time interval $[t_1, t_2),$ where on the contrary only the delay feedback $B_2$ is present.
First, on $[0,t_1],$ since $B_2\equiv 0,$ we are in the situation of \cite{GenniSicon}. Thus, for initial data $u_0\in V$ and $u_1\in H,$ the solution $u$ belongs to $C([0,t_1]; V)\cap C([0,t_1];H).$ Then, we decompose the second interval $[t_1, t_2)$ into the intervals $(t_1+l\tau , t_1+ (l+1)\tau ),$ for $l=0,\dots, L,$ where $L$ is the first value
such that $t_1+(L+1)\tau\ge t_2.$ The last interval is then $(t_1+L\tau , t_2).$
Now, we look at the interval $(t_1, t_1+\tau ).$ In this time interval problem (\ref{1.1}) -- (\ref{1.2}) can be rewritten as
\begin{equation}\label{pbritardo}
\begin{array}{l}
u_{tt}(t)+ Au(t)= g_1(t) +f(u), \quad t\in (t_1, t_1+\tau ),\\
u(t_1+)=u(t_1-) \ \mbox{\rm and }\ u_t(t_1+)=u_t(t_1-),
\end{array}
\end{equation}
where $g_1(t)=-B_2(t) B_2^*(t) u_t(t-\tau)$ belongs to $C([t_1, t_1+\tau); H)$
from the first step. Indeed, for $t\in (t_1, t_1+\tau ),$ it is $t-\tau \in (0, t_1).$
Then, since $(u(t_1-), u_t(t_1-))$ belongs to $V\times H,$  the existence of local  solution  $u\in C^1([t_1, t_1+\delta ]; H)\cap C([t_1, t_1+\delta ]; V)\,,$ $\delta\le\tau,$ follows from \cite[Theorems 1.4 and 1.5, Ch. 6]{Pazy}.
Now observe that, from (\ref{stimaderD2abstrait}),
$$E(t)\le e^{2M_1\tau }E(t_1),\quad \forall \ t\in [t_1, t_1+\delta]\,,$$
then $\delta =\tau\,,$ namely there exists a solution
$u\in C^1([t_1, t_1+\tau ]; H)\cap C([t_1, t_1+\tau ]; V)\,.$
By iterating this procedure we find $u\in  C^1([t_1+\tau, t_1+2\tau]; H)\cap C([t_1+\tau, t_1+2\tau ]; V)$ and then on the whole interval $(t_1, t_2).\qed$

\section{Stability results}
\label{st}

\hspace{5mm}

\setcounter{equation}{0}

In this section, we give sufficient conditions ensuring stability results in case of distributed/lo\-ca\-lized damping.
More explicit conditions are given in the linear case, improving previous results given in \cite{NP2014}.

\subsection{Distributed damping}

\label{DD}

\hspace{5mm}

First of all, consider the case $U_1=W=H,$ that is the case of distributed feedback $B_1(t).$ In  this case, of course, the constant $C$ in the estimate (\ref{embedding}) is $1\,.$
The following result holds (see \cite[Theorem 4.1]{GenniSicon}; cfr. \cite{HMV}).

\begin{Theorem} \label{ObservST}
Assume  {\rm i)} and $(\ref{NL2}).$ Then, any solution $u$ of $(\ref{1.1})-(\ref{1.2})$ satisfies
\begin{equation}\label{obsshort}
E_S(t_{2n+1})\le \frac 1 {1+\frac {T_{2n}^3}{30}\frac 1 {  \frac 4 {\lambda_1 m_{2n}}+\frac {3T_{2n}^2}{32 m_{2n}}+\frac {M_{2n}T_{2n}^2}{16\lambda_1}} }E_S(t_{2n}),\ n\in\nat,
\end{equation}
where $\lambda_1$ is the constant in $(\ref{autoval}).$
\end{Theorem}

\begin{Theorem}\label{First}
Assume {\rm i)}, {\rm ii)}, $(\ref{autoval})$ and $(\ref{NL1})$ -- $(\ref{length}).$
If
\begin{equation}\label{M33}
  \sum_{n=0 }^\infty
(2
M_{2n+1}T_{2n+1}+\ln \tilde c_{n})=-\infty,
\end{equation}
where
\begin{equation}\label{cn}
\tilde c_n=
 \frac 1 {1+\frac {T_{2n}^3}{30}\frac 1 {  \frac 4 {\lambda_1 m_{2n}}+\frac {3T_{2n}^2}{32 m_{2n}}+\frac {M_{2n}T_{2n}^2}{16\lambda_1}} }
  +
  M_{2n+1}T_{2n+1},
  \end{equation}
then
system $(\ref{1.1})-(\ref{1.2})$ is asymptotically stable, that is any solution $u$ of $(\ref{1.1})-(\ref{1.2})$ satisfies $E_S(t)\rightarrow 0$ as $t\rightarrow +\infty.$
\end{Theorem}

\noindent For some comments on  \eqref{M33} we refer to the next Remark \ref{altro}.

\noindent {\bf Proof of Theorem \ref{First}.}
Observe that (\ref{stimaderD2abstrait}) implies
$$E^{\prime}(t)\le 2M_{2n+1}E(t),\quad t\in I_{2n+1}=[t_{2n+1},t_{2n+2}),\ n\in\nat.$$
Then we deduce
\begin{equation}\label{M12A}
E(t_{2n+2})\le e^{2M_{2n+1}T_{2n+1}}E(t_{2n+1}),
\quad \forall \ n\in\nat.
\end{equation}
Now, note that

$$
E(t_{2n+1} ) =\ds{E_S(t_{2n+1})+ \frac{1}{2} \int_{t_{2n+1} -\tau}^{t_{2n+1}} \|B^*_2(s+\tau)u_t(s)\|^2_{U_2}ds},
$$
and then, as $\vert I_{2n}\vert \ge \tau,$ $n\in\nat,$ and $B_2(t)$ is null on the intervals $I_{2n},$
\begin{equation}\label{ven}
\begin{array}{l}
\displaystyle{
E(t_{2n+1}) \le E_S(t_{2n+1})+ \frac{1}{2}M_{2n+1} \int_{t_{2n+1} -\tau}^{\min \{t_{2n+1}, t_{2n+2}-\tau\}} \|u_t(s)\|^2_H ds}\\\medskip
\hspace{1.6 cm}\displaystyle{
\le E_S(t_{2n+1})+ M_{2n+1} \int_{t_{2n+1} -\tau}^{\min \{t_{2n+1}, t_{2n+2}-\tau\}} E_S(s) ds}\\\medskip
\hspace{1.6 cm}\displaystyle{
\le E_S(t_{2n+1})+
M_{2n+1} T_{2n+1}E_S(t_{2n+1}-\tau )}\\\medskip
\hspace{1.6 cm}\displaystyle{
\le E_S(t_{2n+1})+
M_{2n+1}T_{2n+1} E_S(t_{2n}).}
\end{array}
\end{equation}
Then, from Theorem \ref{ObservST} and (\ref{ven}) we deduce
\begin{equation}\label{ven1}
E(t_{2n+1})\le
\left(
 \frac 1 {1+\frac {T_{2n}^3}{30}\frac 1 {  \frac 4 {\lambda_1 m_{2n}}+\frac {3T_{2n}^2}{32 m_{2n}}+\frac {M_{2n}T_{2n}^2}{16\lambda_1}} }
  +
  M_{2n+1}T_{2n+1}
\right) E_S(t_{2n}),
\end{equation}
and therefore, by (\ref{M12A}),
\begin{equation}\label{ven2}
\begin{array}{l}
E_S(t_{2n+2})\le e^{2M_{2n+1}T_{2n+1}}E(t_{2n+1})\\\medskip
\hspace{1.8 cm}
\displaystyle{
\le
e^{2M_{2n+1}T_{2n+1}}
\left(
 \frac 1 {1+\frac {T_{2n}^3}{30}\frac 1 {  \frac 4 {\lambda_1 m_{2n}}+\frac {3T_{2n}^2}{32 m_{2n}}+\frac {M_{2n}T_{2n}^2}{16\lambda_1}} }
  +
 M_{2n+1}T_{2n+1}
\right)
 E_S(t_{2n}).}
\end{array}
\end{equation}
Since (\ref{ven2}) holds for any $n\in\nat$ we conclude
\begin{equation}\label{ven4}
E_S(t_{2n+2})\le\displaystyle{\Pi_{p=0}^n} e^{2M_{2p+1}T_{2p+1}}
\left(
 \frac 1 {1+\frac {T_{2p}^3}{30}\frac 1 {  \frac 4 {\lambda_1 m_{2p}}+\frac {3T_{2p}^2}{32 m_{2p}}+\frac {M_{2p}T_{2p}^2}{16\lambda_1}} }
  +
   M_{2p+1}T_{2p+1}
\right)
 E_S(0).
\end{equation}
Now observe that the standard energy $E_S(\cdot )$ is not decreasing in general. However, it is decreasing
for $t\in [t_{2n}, t_{2n+1})$, when only the standard dissipative damping acts and so
\begin{equation}\label{oggi3}
E_S(t)\le E_S(t_{2n}),\quad \forall \; t\in [t_{2n}, t_{2n+1}).
\end{equation}
Moreover, for $t\in [t_{2n+1}, t_{2n+2}),$ it results
\begin{equation}\label{oggi4}
E_S(t)\le E(t)\le e^{2M_{2n+1}T_{2n+1}}E(t_{2n+1}),
\end{equation}
where in the second inequality we have used  (\ref{stimaderD2abstrait}).

Then, by (\ref{ven4}), (\ref{oggi3}), (\ref{oggi4}) and (\ref{ven1}), asymptotic stability occurs  if
(\ref{M33}) is satisfied.\qed

\begin{Remark}\label{altro}
{\rm Observe that (\ref{M33}) holds true if the following easier conditions are satisfied:
 \begin{equation}\label{prima}
\sum_{n=0}^\infty M_{2n+1}T_{2n+1}<+\infty
\end{equation}
and
\begin{equation}\label{seconda}
\sum_{n=0}^\infty
\ln
\left(
1+\frac {T_{2n}^3}{30}
\frac 1
{  \frac 4 {\lambda_1 m_{2n}}+\frac {3T_{2n}^2}{32 m_{2n}}+\frac {M_{2n}T_{2n}^2}{16\lambda_1} }
\right )=+\infty.
\end{equation}
Indeed,  it is easy to see (cfr. \cite{Pignotti15}) that (\ref{prima}) and
\begin{equation}\label{terza}
\sum_{n=0}^\infty \ln \tilde c_n =-\infty
\end{equation}
with $\tilde c_n,$ $n\in\nat,$ as in (\ref{cn}), imply (\ref{M33}).
Now it is sufficient to observe that, under assumption (\ref{prima}), the conditions (\ref{seconda}) and (\ref{terza}) are equivalent.
Indeed if (\ref{terza}) holds true
then
$$-\ln
\left(
1+\frac {T_{2n}^3}{30}\frac 1 {  \frac 4 {\lambda_1 m_{2n}}+\frac {3T_{2n}^2}{32 m_{2n}}+\frac {M_{2n}T_{2n}^2}{16\lambda_1} }
\right )=
 \ln \left(
\frac 1
{1+\frac {T_{2n}^3}{30}\frac 1 {  \frac 4 {\lambda_1 m_{2n}}+\frac {3T_{2n}^2}{32 m_{2n}}+\frac {M_{2n}T_{2n}^2}{16\lambda_1} }}
\right)< \ln \tilde c_n\,$$
and therefore also (\ref{seconda}) is satisfied. Assume now that (\ref{prima}) and (\ref{seconda}) are satisfied. Then, by (\ref{prima}),
\begin{equation}\label{nec}
M_{2n+1}T_{2n+1}\rightarrow 0,\quad n\rightarrow \infty.
\end{equation}
If {(\ref{terza})} does not hold then it has to be
$$\ln \tilde c_n \; \rightarrow 0,
\quad n\rightarrow \infty.$$
But then, by (\ref{nec}), it results
$$
\tilde c_n \sim \frac 1 {1+\frac {T_{2n}^3}{30}\frac 1 {  \frac 4 {\lambda_1 m_{2n}}+\frac {3T_{2n}^2}{32 m_{2n}}+\frac {M_{2n}T_{2n}^2}{16\lambda_1}} },$$
in contradiction with {(\ref{seconda})}.
}
\end{Remark}

\begin{Remark}\label{altraenergia}{\rm
Observe that, under the assumptions of Theorem \ref{First}, one can prove that
also
\begin{equation}\label{asymptoticenergy}
E(t) \rightarrow 0 \quad \mbox{\rm as } t\rightarrow +\infty,
\end{equation}
for every solution $u$ of $(\ref{1.1})-(\ref{1.2}).$
Indeed,
recall that
$$E(t)= E_S(t)+\frac {1} 2 \int_{t-\tau }^t\Vert B_2^*(s+\tau )u_t(s)\Vert_{U_2}^2 ds$$
  and that we are assuming $T_{2n}\ge\tau,$ for all $n\in\nat,$ and  $B_2$  null in the intervals $I_{2n}.$
  Then, if $t\in \bar I_{2n}=[t_{2n}, t_{2n+1}],$
  \begin{equation}\label{utile}
  \begin{array}{l}
  \displaystyle{
  E(t)\le E_S(t)+\frac {1} 2 \int_{t_{2n+1}-\tau }^{\min \{ t, t_{2n+2}-\tau\}}\Vert B_2^*(s+\tau )u_t(s)\Vert_{U_2}^2 ds}\\
  \displaystyle{
  \hspace{1 cm}\le E_S(t)+  M_{2n+1}\int_{t_{2n+1}-\tau }^{\min \{ t, t_{2n+2}-\tau\}} E_S(s) ds
 \le E_S(t)+   M_{2n+1}T_{2n+1}E_S(t_{2n});}
 \end{array}
 \end{equation}
if $t\in I_{2n+1},$ by using (\ref{stimaderD2abstrait}) and (\ref{utile}), we have
  \begin{equation}\label{utile2}
  E(t)\le e^{2M_{2n+1}T_{2n+1}}
  E(t_{2n+1})
  \le e^{2M_{2n+1}T_{2n+1}}[E_S(t_{2n+1})+  M_{2n+1}T_{2n+1}E_S(t_{2n})]\,.
  \end{equation}
 Therefore, observing that by (\ref{M33}) one has
 $$\sup_n M_{2n+1}T_{2n+1}<+\infty,$$
 $(\ref{asymptoticenergy})$ follows when Theorem \ref{First} applies.
 }\end{Remark}

\vspace{0.5cm}

 \noindent As an example of application of Theorem \ref{First} one can consider problem (\ref{Wo.1}) -- (\ref{Wo.3}) assuming

$i_w$)  $0< m_{2n}\le  b_1(t)\le M_{2n}$,
$b_2(t)=0,$
for all $t\in I_{2n}=[t_{2n},t_{2n+1})$ and $b_1\in C^1(\bar I_{2n})$ for all $\ n\in\nat;$

\vskip+5 pt

 $ii_w$) $\vert b_2(t)\vert \le M_{2n+1}$,
  $b_1(t)=0$
for all $t\in I_{2n+1}=[t_{2n+1},t_{2n+2})$ and $b_2\in C^1(\bar I_{2n+1})$ for all  $\ n\in\nat.$

\vspace{0.5cm}
The previous result can be extended to a more general situation. Indeed, consider
the nonlinear wave system
 \begin{eqnarray}
& &u_{tt}(t) +A u (t)+B_1(t)B_1^*(t)g(u_t)+B_2(t)B_2^*(t)u_t(t-\tau) =f(u),\quad t>0,\label{g1}\\
& &u(0)=u_0\quad \mbox{\rm and}\quad u_t(0)=u_1,\quad\label{g2}
\end{eqnarray}
with $(u_0, u_1) \in V\times H$.
On the functions $g$ and $f$
 we make the following assumptions:
\[
(A)\quad \left\{\begin{array}{l}
g:\RR\longrightarrow \RR\mbox{ is a }C^1 \mbox{ function with $g(0)=0$},\\
\exists\,B\geq A>0\mbox{ such that $0<A\leq g'(v)\leq B$
$\forall\,v\in\RR$},\\
\mbox{$f$ satisfies \eqref{NL1} and \eqref{NL2}}.\end{array}\right.
\]
Moreover, on $B_1$ we assume, in place of i), that, for all $n\in\nat$, there exist positive constants $m_{2n}$,
 $M_{2n},$  with $m_{2n}\leq M_{2n}$, such that, for all $u\in H$, we have

\vskip+5pt

i') $m_{2n}\|u\|_W^2\le  \langle B_1^*(t)u, B_1^*(t) g(u)\rangle_{U_1}\le M_{2n} \|u\|_W^2$ for $t\in I_{2n}=[t_{2n},t_{2n+1}),$  $\ \forall\ n\in\nat.$

\vspace{0.5cm}

As prototype, one can think to the problem
 \begin{eqnarray}
&& u_{tt}(x,t) -\Delta u (x,t)+b_1(t) g(u_t)+b_2(t)  u_t(x,t-\tau)= -|u|^pu \ \mbox{\rm in}\  \Omega \times (0,+\infty),\\
&&u(x,t)=0 \ \mbox{\rm in}\  \partial\Omega\times  (0,+\infty),\\
& &u(x,0)=u_0(x)\quad \mbox{\rm and}\quad u_t(x,0)=u_1(x)\quad \hbox{\rm
in} \quad \Omega,
 \end{eqnarray}
 where $\Omega$,  $u_0$, $u_1$, $b_1$, $b_2$ and $p$ are as before.

For \eqref{g1} -- \eqref{g2}, if $B_2(t)=0$ for all $t \in (0,+\infty)$, Theorem \ref{ObservST} becomes (see \cite[Theorem 5.1]{GenniSicon})

\begin{Theorem} \label{ObservSTgenerale}
Assume {\rm i') } and suppose that also $(A)$ holds. Then, any solution $u$ of $(\ref{g1})-(\ref{g2})$ satisfies
\begin{equation}\label{obsshort_1}
E_S(t_{2n+1})\le \frac 1 {1+\frac {T_{2n}^3}{30}\frac 1 {  \frac 4 {\lambda_1 m_{2n}}+\frac {3T_{2n}^2}{32 m_{2n}}+\frac {M_{2n}T_{2n}^2}{16\lambda_1}} }E_S(t_{2n}),\ n\in\nat.
\end{equation}
\end{Theorem}

Therefore, observe that, since $B_1(t)=0$ for all $t \in I_{2n+1}$, \eqref{stimaderD2abstrait} still holds. Thus, as Theorem \ref{ObservST} implies Theorem \ref{First},  Theorem \ref{ObservSTgenerale} immediately gives, for global defined solutions, the following
fundamental application via \eqref{stimaderD2abstrait}:

\begin{Theorem}\label{First_generale}
Assume {\rm  i')} and {\rm ii)}. Moreover suppose that
$(\ref{length})$ and $(A)$ are satisfied.
If \eqref{M33} holds, then
system $(\ref{g1})-(\ref{g2})$ is asymptotically stable, that is any  solution $u$ of $(\ref{g1})-(\ref{g2})$ satisfies $E_S(t)\rightarrow 0$ as $t\rightarrow +\infty.$
\end{Theorem}

\begin{Remark}\label{rem_generale}
{\rm
Observe that problems  \eqref{1.1} --  \eqref{1.2} and $(\ref{g1})-(\ref{g2})$ with  $B_2(t)\equiv 0$ correspond to the models with on-off damping considered in \cite{GenniSicon}. Hence Theorems \ref{First} and \ref{First_generale} give stability results also in these situations.
}\end{Remark}

Of course, the abstract setting of the previous theorems let us deal with higher order problems in bounded and smooth domain of $\RR^N$. For example, Theorem \ref{First_generale} can be applied to the problem
\begin{eqnarray}
&& u_{tt}(x,t) +\Delta ^{2m}u (x,t)+b_1(t) g(u_t)+b_2(t)  u_t(x,t-\tau)= f(u) \ \mbox{\rm in}\  \Omega \times (0,+\infty),\\
&&Cu(x,t)=0 \in \RR^{2m} \ \mbox{\rm in}\  \partial\Omega\times  (0,+\infty),\\
& &u(x,0)=u_0(x) \in D(\Delta^m)\quad \mbox{\rm and}\quad u_t(x,0)=u_1(x)\quad \hbox{\rm
in} \quad \Omega,
 \end{eqnarray}
 where $m \in N$, $f$, $g$, $b_1$, $b_2$ and $p$ are as before and $C$ is a boundary operator such that {\it the first eigenvalue of $\Delta^{2m}$ under the boundary conditions $Cu(x,t)=0 \in \RR^{2m}$ in $\partial \Omega \times (0, +\infty)$ is strictly positive.} For example, one can consider as $C$ the Dirichlet operator, while the case of Neumann boundary conditions must be excluded since the first eigenvalue is $0$.

\subsection{Localized damping}

\label{LD}

\hspace{5mm}

In this section we consider the more general situation $U_1\ne W$. In practice, for concrete models, the feedback operators $B_1$
and $B_2$ may be localized in  subregions of $\Omega.$

\begin{Proposition}\label{derivECris}
Assume ${\rm  i)}$, ${\rm ii)},$ $(\ref{NL1})$ and $(\ref{length}).$
For any regular solution of problem $(\ref{1.1})-(\ref{1.2})$ the energy $E_S$ is decreasing
on the intervals $I_{2n},$ $n\in\nat.$ In particular,
\begin{equation}\label{stimader2Cris3}
E_S^{\prime}(t)= -
            \|B_1^*(t)u_t(t)\|^2_{U_1}.
\end{equation}
Moreover, on the intervals $I_{2n+1},$ $n\in\nat,$
the estimate $(\ref{stimaderD2abstrait})$ holds.
\end{Proposition}

\noindent{\bf Proof.} By differentiating $E_S(\cdot ),$ we have

$$E_S^\prime (t)=\langle u_t, u\rangle_V+\langle u_{tt}, u_t\rangle_H-\langle f(u), u_t\rangle_H.$$
Then, recalling that $B_2(t)=0$ in $I_{2n}$, from equation $(\ref{1.1})$ it follows that
$$E_S^\prime (t)=\langle u_t, u_{tt}+Au-f(u)\rangle_{V,V^\prime }=-\langle u_t, B_1(t) B_1^*(t) u_t(t)\rangle_{V,V^\prime, }$$
for all $t \in I_{2n}$.
Thus, identity $(\ref{stimader2Cris3})$ holds.\qed

Consider now the system
\begin{eqnarray}
& &w_{tt}(t) +A w (t)+B_1(t) B_1^*(t) w_t=f(w),\quad t\in (t_{2n}, t_{2n+1}),\ n\in\nat, \label{cons1abstraitU}\\
& &w(t_{2n})=w_0^n\quad \mbox{\rm and}\quad w_t(t_{2n})=w_1^n\quad\label{cons2NL}
\end{eqnarray}
with $(w_0^n,w_1^n)\in V\times H$.
For our stability result we need that the next    observability type  inequality holds.
Namely we assume that, for every $n$,  there exists a time $\overline T_n$
such that
\begin{equation}\label{T2nNL}
T_{2n}>\overline T_n,
\end{equation}
and that,
for every $n$ and every time $T,$ with $T_{2n}\ge T> \overline T_n,$
 there is     a constant $d_n,$ depending on $T$ but independent of $(w_0^n,w_1^n),$ such that
\begin{equation}\label{OCabstraitCrisNL}
E_S(t_{2n}+T)\le d_n\int_{t_{2n}}^{t_{2n}+T}\Vert B_1^*(t) w_t(t)\Vert_{U_1}^2  dt,
\end{equation}
for every weak solution of problem  \eqref{cons1abstraitU} -- \eqref{cons2NL}
with initial data $(w_0^n,w_1^n)\in V\times H.$

\begin{Remark}{\rm The observability inequality above is satisfied for solutions of wave--type equations when the nonlinearity $f$ satisfies some requirements. For instance in \cite{Zuazua} Zuazua proved (\ref{OCabstraitCrisNL}) if $f$ is globally Lipschitz, as a perturbation of the well--known linear case, or also when $f$ satisfies
\begin{equation}\label{zuazua}
(2+\delta ) F(s)\ge s f(s),
\end{equation}
for some $\delta >0\,$.}
\end{Remark}

\begin{Proposition}\label{mista}
Assume $\mbox{\rm i)}.$
Moreover, we assume that there is a sequence $\{\overline T_n\}_n,$ such that
$(\ref{T2nNL})$ is satisfied and the
observability estimate $(\ref{OCabstraitCrisNL})$ holds for every $T\in (\overline T_n, T_{2n}],$ $\forall\ n\in\nat.$
Then, for any solution of system $(\ref{1.1})-(\ref{1.2})$ we have
\begin{equation}\label{cris1NL}
E_S(t_{2n+1})\le \hat d_n E_S(t_{2n}),\quad \forall\ n\in\nat,
\end{equation}
where
\begin{equation}\label{nnewU}
\displaystyle{
\hat d_n=
\frac{d_n}{d_{n}+1}
},
\end{equation}
$d_n$ being the observability constant in
$(\ref{OCabstraitCrisNL})$ corresponding to the time $T_{2n}.$
 \end{Proposition}

\noindent {\bf Proof.} To prove (\ref {cris1NL}) it suffices to
use the estimate $(\ref{stimader2Cris3})$ in
$(\ref{OCabstraitCrisNL}),$ reminding that $B_2(t)=0$ on $(t_{2n},t_{2n+1})$.
Indeed, $(\ref{OCabstraitCrisNL})$ gives
\begin{equation}\label{Cris4}
E_S(t_{2n+1})\le d_n\int_{t_{2n}}^{t_{2n+1}} \Vert B_1^*(t) u_t(t)\Vert_{U_1}^2 dt.
\end{equation}
By integrating $(\ref{stimader2Cris3})$ on the interval $[t_{2n}, t_{2n+1}],$ we have
$$E_S(t_{2n+1})-E_S(t_{2n})=-\int _{t_{2n}}^{t_{2n+1}} \Vert B_1^*(t) u_t(t)\Vert_{U_1}^2 dt,$$
and therefore, using (\ref {Cris4}),
$$E_S(t_{2n+1})-E_S(t_{2n})\le -\frac 1 {d_n} E_S(t_{2n+1}).$$
Thus,
$$E_S(t_{2n+1})\Big (\frac {d_n+1} {d_n}\Big )\le E_S(t_{2n}).\quad\qed$$

\begin{Theorem}\label{stab2Cris5}
Assume hypotheses of Proposition  $\ref{mista},$ ${\rm ii)},$ $(\ref{NL1})$  and $(\ref{length}).$ If
\begin{equation}\label{starstarNL}
\sum_{n=0}^\infty  [2CM_{2n+1}T_{2n+1}+\ln (\hat d_n+CM_{2n+1}T_{2n+1})]=-\infty,
\end{equation}
where $C$ is the constant in the norm embedding $(\ref{embedding})\,,$
then
system $(\ref{1.1})-(\ref{1.2})$ is asymptotically stable, that is any solution  $u$ of $(\ref{1.1})-(\ref{1.2})$ satisfies
$E_S(t)\rightarrow 0$ as $t\rightarrow +\infty.$
\end{Theorem}

\noindent {\bf Proof.} The proof is analogous to the one of Theorem \ref{First}. Simply we use now  the inequality (\ref{cris1NL}) in place of
(\ref{obsshort}) on the intervals $I_{2n}$.\qed

\begin{Remark}{\rm As in Remark \ref{altro}, one can show
 that (\ref{starstarNL}) is verified if, in particular,
 \begin{equation}\label{true}
\sum_{n=0}^\infty M_{2n+1}T_{2n+1} <+\infty\quad \mbox{\rm and}\ \quad
\sum_{n=0}^\infty \ln \hat d_n=-\infty.
\end{equation}
Observe also that $d_n$ depends on  $n$ since, by hypothesis,  $B_1$ may depend on the time variable. However,
if $B_1$ is independed of $t,$ then by a translation of $t_{2n}$
the constant   $d_n$ becomes independent of $n$. But, if $d_n=d>0$ for all $n$, then the condition
$$
\sum_{n=0}^\infty\ln  \hat d_n =-\infty
$$
is clearly satisfied. On the other hand, the first condition in (\ref{true}) depends only on the length of the intervals $I_{2n+1}$ and on the boundedness constant of $B_2^*$ on the same intervals,
hence (\ref{true}) can be easily checked.
}
\end{Remark}

As an example of model for which this result holds, we can consider
 \begin{eqnarray}
&&u_{tt}(x,t) -\Delta u (x,t)+b_1(t)\chi_\omega u_t(x,t)+b_2(t) \chi_{\tilde\omega} u_t(x,t-\tau)=f(u)\quad \mbox{\rm in}\ \Omega\times
(0,+\infty),\;\;\;\;\;\;\;\;\label{Wo.12}\\
&&u (x,t) =0\quad \mbox{\rm on}\quad\partial\Omega\times
(0,+\infty),\label{Wo.22}\\
&&u(x,0)=u_0(x)\quad \mbox{\rm and}\quad u_t(x,0)=u_1(x)\quad \hbox{\rm
in}\quad\Omega,\label{Wo.32}
\end{eqnarray}
with
initial
data $(u_0, u_1)\in H^1_0(\Omega)\times L^2(\Omega),$
$b_1,b_2$ as before
and the nonlinearity $f$ as in \cite{Zuazua}. Moreover, we assume that the set $\omega \subset \Omega$ satisfies a control geometric property (see \cite{BLR}) and that $\tilde\omega \subset \omega.$

\vspace{0.5cm}
As for the distributed damping, the previous result can be extended to a more general situation. Indeed, consider again
the nonlinear wave system \eqref{g1} -- \eqref{g2}
with $(u_0, u_1) \in V\times H$.
On the functions $g$ and $f$  we assume $(A)$.

 Proposition \ref{derivECris} becomes
 \begin{Proposition}\label{derivECris_1}
Assume $\mbox{\rm i'),\ ii)},$  $(\ref{NL1})$ and $(\ref{length}).$
For any regular solution of problem $(\ref{g1})-(\ref{g2})$ the energy $E_S$ is such that
\begin{equation}\label{stimader2Cris3_1}
E_S^{\prime}(t)= -
            \langle B_1^*(t)u_t(t), B_1^*(t) g(u_t) \rangle_{U_1},
\end{equation}
for all $t \in I_{2n},$ $n\in\nat$.
Moreover, on the intervals $I_{2n+1},$ $n\in\nat,$
the estimate $(\ref{stimaderD2abstrait})$ holds.
\end{Proposition}

As before, consider now the system
\begin{eqnarray}
& &w_{tt}(t) +A w (t)+B_1(t) B_1^*(t) g(w_t)=f(w),\quad t\in (t_{2n}, t_{2n+1}),\ n\in\nat, \label{cons1abstraitU_1}\\
& &w(t_{2n})=w_0^n\quad \mbox{\rm and}\quad w_t(t_{2n})=w_1^n,\quad\label{cons2NL_1}
\end{eqnarray}
with $(w_0^n,w_1^n)\in V\times H.$
For our stability result we need that the next  inequality holds.
Namely we assume that, for every $n$, there exists a time ${\overline T}_n$ such that \eqref{T2nNL} holds
and that,
for every $n$ and every time $T,$ with $T_{2n}\ge T> \overline T_n,$
 there is     a constant $d_n,$ depending on $T$ but independent of $(w_0^n,w_1^n),$ such that
\begin{equation}\label{OCabstraitCrisNL_1}
E_S(t_{2n}+T)\le d_n E_S(t_{2n}),
\end{equation}
for every weak solution of problem  \eqref{cons1abstraitU_1} -- \eqref{cons2NL_1}
with initial data $(w_0^n,w_1^n)\in V\times H.$

The  inequality above is satisfied for solutions of wave--type equations when the
nonlinearities $f$ and $g$ satisfy some requirements. In \cite{Zuazua}, for example, (\ref{OCabstraitCrisNL_1}) is proved if $f$ is globally Lipschitz or when $f$ satisfies
\eqref{zuazua}
for some $\delta >0$ and
$g$ is globally Lipschitz (hence if $g$ is as in $(A)$) and there exists $c>0$ such that
\[
g(s)s \ge c|s|^2, \quad \forall s \in \RR.
\]

Theorem \ref{stab2Cris5} becomes
\begin{Theorem}\label{stab2Cris5_1}
Assume $\mbox{\rm i'), ii)},$ $(\ref{NL1})$ and $(\ref{length}).$
Moreover, we assume that there is a sequence $\{\overline T_n\}_n,$ such that
$(\ref{T2nNL})$ is satisfied and  $(\ref{OCabstraitCrisNL_1})$ holds for every $T\in (\overline T_n, T_{2n}],$ $\forall\ n\in\nat.$
If \eqref{starstarNL} holds
then
system $(\ref{g1})-(\ref{g2})$ is asymptotically stable, that is any solution  $u$ of $(\ref{g1})-(\ref{g2})$ satisfies
$E_S(t)\rightarrow 0$ as $t\rightarrow +\infty.$
\end{Theorem}

\begin{Remark}{\rm
One can make the same considerations made in Remark \ref{rem_generale}, obtaining stability results  also for the localized on-off damping. These results are then more general than the ones proved in \cite{GenniSicon}.
}\end{Remark}

\subsection{Localized damping: the linear case}

\label{LDLinear}

\hspace{5mm}

In the linear case (i.e. $f\equiv 0$) we can improve previous results given in \cite{NP2014} by removing the assumption (\ref{new}) on the coefficients.
As in \cite{NP2014} we can determine more explicitely, in terms of the coefficients $T_{2n}, m_{2n}, M_{2n},$ the constant $\hat d_n$ of Proposition \ref{mista}, for all $n\in\nat\,.$

Consider now the conservative system associated with (\ref{1.1}) -- (\ref{1.2})
 \begin{eqnarray}
& &w_{tt}(t) +A w (t)=0,\quad t>0,\label{cons1}\\
& &w(0)=w_0\quad \mbox{\rm and}\quad w_t(0)=w_1,\quad\label{cons2}
\end{eqnarray}
with $(w_0,w_1)\in V\times H.$

To prove stability results  we need that a suitable   observability inequality holds.
Then, we assume that there exists a time $\overline T>0$ such that, for every time $T>\overline T$, there is     a constant $c,$ depending on $T$ but independent of the initial data, such that
\begin{equation}\label{OC}
E_S(0)\le c\int_0^{T}\|w_t(s)\|^2_W  ds,
\end{equation}
for every weak solution of problem $(\ref{cons1})-(\ref{cons2})$
with initial data $(w_0,w_1)\in V\times H.$

The following result is  proved in \cite{NP2014}:

\begin{Proposition}\label{obsexplicit}
Assume $\mbox{\rm i)}$ and $f\equiv 0\,.$
Moreover,
we assume that the observability inequality $(\ref{OC})$ holds for every  time $T>\overline T$
and that, setting  $T^*:=\inf_n \{T_{2n}\},$
\begin{equation}\label{ampiezze}
T^*>\overline T.
\end{equation}
Then, for any solution of system $(\ref{1.1})-(\ref{1.2})$ we have
\begin{equation}\label{stimabuona2Anew}
E_S(t_{2n+1})\le \hat c_n E_S(t_{2n}),\quad \forall\ n\in\nat,
\end{equation}
where
\begin{equation}\label{nnew}
\displaystyle{
\hat c_n=
\frac{2c(1+4C^2T_{2n}^2M_{2n}^2) }
{m_{2n}+ 2c(1+4C^2T_{2n}^2M_{2n}^2)}
},
\end{equation}
  $c$ being the observability constant in $(\ref{OC})$ corresponding to the time $T^*$ and $C$ the constant in the norm embedding $(\ref{embedding})$  between $W$ and $H$.
\end{Proposition}

Combining the previous proposition with estimate (\ref{stimaderD2abstrait}) one can obtain the following theorem.

\begin{Theorem}\label{linear}
Assume hypotheses of Proposition $\ref{obsexplicit},$
 ${\rm ii)},$ $(\ref{NL1})$ and $(\ref{length}).$
If
\begin{equation}\label{luglio1}
\sum_{n=0}^\infty [2 CM_{2n+1}T_{2n+1}+\ln (\hat c_n +CM_{2n+1}T_{2n+1})]=-\infty,
\end{equation}
where  $\hat c_n$ is as in $(\ref{nnew})$ and $C$ is the constant in the norm embedding $(\ref{embedding}),$
then system $(\ref{1.1})-(\ref{1.2})$ is asymptotically stable, that is for every solution of $(\ref{1.1})-(\ref{1.2})$ $E_S(t)\rightarrow 0$ as $t\rightarrow +\infty\,.$
\end{Theorem}

\noindent {\bf Proof.} The proof is analogous to the one of Theorem \ref{First}. Simply we use now the inequality (\ref{stimabuona2Anew}) in place  of
(\ref{obsshort}) on the intervals $I_{2n}$.\qed

\begin{Remark}\label{remarklinear}{\rm
As in Remark \ref{altro} we can show that (\ref{luglio1}) is verified in particular
if
\begin{equation}\label{luglio2}
\sum_{n=0}^\infty M_{2n+1}T_{2n+1}<+\infty ,\quad \mbox{\rm and }\quad \sum_{n=0}^\infty \ln\hat c_n=-\infty.
\end{equation}
Now, it is easy to see that the second condition of (\ref{luglio2}) is equivalent  (see the proof of \cite[Theorem 3.3]{NP2014} for details)
to
\begin{equation}\label{luglio3}
\sum_{n=0}^\infty \frac {m_{2n}}
{1+4C^2T^2_{2n}M^2_{2n}}=+\infty.
\end{equation}
which is, together with the first condition of (\ref{luglio2}) on the intervals with delay, the assumption of \cite[Theorem 3.3]{NP2014}. Actually, as clearly appears from the proof, Theorem 3.3 of \cite{NP2014} holds true under the more general condition (\ref{luglio1}). The authors there preferred, for sake of clairness, to formulate the assumption in an easier but less general form.}
\end{Remark}

\begin{Remark}{\rm
Observe that Theorem \ref{linear}  significantly improve \cite[Theorem 3.3]{NP2014}. Indeed it allows to obtain the same stability result by removing the assumption (\ref{new}) on the coefficients, which is crucial in the proof of  \cite[Theorem 3.3]{NP2014}.
}\end{Remark}

\section{Stability result:  localized positive--negative damping without delay}

\label{PN}

\setcounter{equation}{0}

In this section we want to generalize the results given in \cite{Genni2} to the localized situation. In particular, in order to deal with a positive--negative damping, we consider the problem
 \begin{eqnarray}
& &u_{tt}(t) +A u (t)+B_1(t)B_1^*(t)u_t(t) - B_3(t)B_3^*(t)u_t(t)=f(u),\quad t>0,\label{4.1}\\
& &u(0)=u_0\quad \mbox{\rm and}\quad u_t(0)=u_1,\quad\label{4.2}
\end{eqnarray}
where
$B_i(t)\in {\mathcal L} (U_i,  {H}),$ $i=1,3.$ Here  $H$ and $U_i,\;i=1, 3,$ are  real Hilbert spaces as before.
On the time--dependent operators $B_i$ we assume
$$B_1^*(t)B_3^*(t)=0,\quad \forall\; t>0$$
and, for all $n\in\nat$, there exists $t_n>0$, with $t_n<t_{n+1}$,
such that
\begin{eqnarray*}
B_3(t)=0,\;  \forall\ t\in I_{2n}=[t_{2n},t_{2n+1}),\\
B_1(t)=0,
\; \forall \  t\in I_{2n+1}=[t_{2n+1},t_{2n+2}),
\end{eqnarray*}
with $B_1\in C^1([t_{2n},t_{2n+1}]; {\mathcal L} (U_1,  {H}))$
and $B_3\in C([t_{2n+1},t_{2n+2}]; {\mathcal L} (U_3,  {H}))$.
We further assume that there exist two Hilbert spaces $W_1, W_3$ such that,
for $i=1,3,$
\begin{equation}\label{embeddi}
 \|u\|_{W_i}^2\le C_i \|u\|_H^2,\quad\forall \;u\in H \ \mbox{\rm with}\ \  C_i>0\ \
\mbox{\rm  independent of}\  u,
\end{equation}
and, for all $n\in\nat$, there exist three positive constants $m_{2n}$,
 $M_{2n}$ and $M_{2n+1}$, with $m_{2n}\leq M_{2n}$, such that for all $u\in H$ we have

j) $m_{2n}\|u\|_{W_1}^2\le  \|B_1^*(t)u\|_{U_1}^2\le M_{2n} \|u\|_{W_1}^2$ for $t\in I_{2n}=[t_{2n},t_{2n+1}),$  $\ \forall\ n\in\nat;$

jj)$\|B_3^*(t)u\|_{U_3}^2 \le M_{2n+1}\|u\|_{W_3}^2 $ for $t\in I_{2n+1}=[t_{2n+1},t_{2n+2}),$ $\ \forall\ n\in\nat.$

\vspace{0.5cm}

The energy functional $E(t)$ coincide in this case with $E_S(t)$ and the next result holds.

\begin{Proposition}\label{derivECris_senzaritardo}
Assume $(\ref{NL1}).$ Then,
for any regular solution of problem $(\ref{4.1})-(\ref{4.2})$ the energy is decreasing
on the intervals $I_{2n}$ and increasing on $I_{2n+1}$, $n\in\nat.$ In particular,
\begin{equation}\label{stimader2Cris3_senzaritardo}
E_S^{\prime}(t)= -
            \|B_1^*(t)u_t(t)\|^2_{U_1}, \quad \forall \; t \in I_{2n}
\end{equation}
and
\begin{equation}\label{stimader2Cris3_senzaritardo1}
E_S^{\prime}(t)=
            \|B_3^*(t)u_t(t)\|^2_{U_3}, \quad \forall\; t \in I_{2n+1}
\end{equation}
\end{Proposition}

\noindent{\bf Proof.} Proceeding as in Proposition \ref{derivECris}, one has

$$E_S^\prime (t)=\langle u_t, u\rangle_V+\langle u_{tt}, u_t\rangle_H-\langle f(u), u_t\rangle_H.$$
Then, recalling that $B_3(t)=0$ in $I_{2n}$ and $B_1(t)=0$ in $I_{2n+1}$, from equation $(\ref{4.1})$ it follows that
$$E_S^\prime (t)=\langle u_t, u_{tt}+Au-f(u)\rangle_{V,V^\prime }=-\langle u_t, B_1(t) B_1^*(t) u_t(t)\rangle_{V,V^\prime }$$
for all $t \in I_{2n}$ and
\[
E_S^\prime (t)=\langle u_t, u_{tt}+Au-f(u)\rangle_{V,V^\prime }=\langle u_t, B_3(t) B_3^*(t) u_t(t)\rangle_{V,V^\prime }
\]
for all $t \in I_{2n+1}$.
Thus, identities $(\ref{stimader2Cris3_senzaritardo})$  and $(\ref{stimader2Cris3_senzaritardo1})$ hold.\qed

As in the previous section we consider the system \eqref{cons1abstraitU} -- \eqref{cons2NL} for which we assume that the observability inequality \eqref{OCabstraitCrisNL}
holds.

Setting again $T_n := t_{n+1}-t_n$, we have:
 \begin{Proposition}\label{mista_senzaritardo}
Assume ${\rm i)}$ and suppose that
there is a sequence $\{\overline T_n\}_n,$ such that
$(\ref{T2nNL})$ and $(\ref{OCabstraitCrisNL})$ hold for every $T\in (\overline T_n, T_{2n}],$ $\forall\ n\in\nat.$
Then, for any solution of system $(\ref{1.1})-(\ref{1.2})$ we have
\begin{equation}\label{cris1NL_1}
E_S(t_{2n+1})\le \hat d_n E_S(t_{2n}),\quad \forall\ n\in\nat,
\end{equation}
where
\begin{equation}\label{nnewU_2}
\displaystyle{
\hat d_n=
\frac{d_n}{d_{n}+1}
},
\end{equation}
$d_n$ being the observability constant in
$(\ref{OCabstraitCrisNL})$ corresponding to the time $T_{2n}.$
 \end{Proposition}

The proof of the previous Proposition is similar to the one of Proposition \ref{mista}, so we omit it.

Using Proposition $\ref{mista_senzaritardo},$ one can give an asymptotic stability result.

\begin{Theorem}\label{stab2Cris5_senzaritardo}
Assume hypotheses of Proposition  $\ref{mista_senzaritardo},$
${\rm ii)},$ $(\ref{NL1})$ and $(\ref{length}).$
If
\begin{equation}\label{starstarNL_1}
\sum_{n=0}^\infty  [2C_3M_{2n+1}T_{2n+1}+\ln \hat d_n ]=-\infty,
\end{equation}
where $C_3$ is the constant in the norm embedding $(\ref{embeddi})$ between $W_3$ and $H,$
then
system $(\ref{4.1})-(\ref{4.2})$ is asymptotically stable, that is any solution  $u$ of $(\ref{4.1})-(\ref{4.2})$ satisfies
$E_S(t)\rightarrow 0$ as $t\rightarrow +\infty.$
\end{Theorem}

\noindent{\bf Proof.} From (\ref{stimader2Cris3_senzaritardo1}) and (\ref{embeddi}) we obtain
$$E_S(t)\le e^{2C_3M_{2n+1}T_{2n+1}} E_S(t_{2n+1}),\quad \forall\ t\in I_{2n+1}=[t_{2n+2}, t_{2n+1}].$$
Therefore, by using $(\ref{cris1NL_1}),$
we have
$$ E_S(t_{2n+2})\le e^{2C_3M_{2n+1}T_{2n+1}} \hat d_n E_S (t_{2n}).$$

Now, we can conclude proceeding as in the proof of Theorem \ref{First}.\qed

\begin{Remark}{\rm
In particular (\ref{starstarNL_1}) is satisfied if
$$
\sum_{n=0}^\infty M_{2n+1}T_{2n+1} <+\infty\quad \mbox{\rm and}\ \quad
\sum_{n=0}^\infty \ln \hat d_n=-\infty.
$$}
\end{Remark}

As an example of model for which the previous result holds, we can consider
 $$
 \begin{array}{l}
u_{tt}(x,t) -\Delta u (x,t)+b_1(t)\chi_\omega u_t(x,t)-b_3(t) \chi_{\tilde\omega} u_t(x,t)=f(u)\quad
\mbox{\rm in}\ \Omega\times
(0,+\infty),\\
u (x,t) =0\quad \mbox{\rm on}\quad\partial\Omega\times
(0,+\infty),\\
u(x,0)=u_0(x)\quad \mbox{\rm and}\quad u_t(x,0)=u_1(x)\quad \hbox{\rm
in}\quad\Omega,
\end{array}
$$
with
initial
data $(u_0, u_1)\in H^1_0(\Omega)\times L^2(\Omega),$
$b_1,b_3$ in $L^\infty (0,+\infty)$ such that
$$b_1(t)b_3(t)=0,\quad\forall\ t>0,$$
and the nonlinearity $f$ as in  \cite{Zuazua}. Moreover, we assume that the set $\omega \subset \Omega$ satisfies a control geometric property.

 On the coefficients $b_1$ and $b_3$ we assume
\vspace{0.5cm}

$j_w$)  $0< m_{2n}\le  b_1(t)\le M_{2n}$,
$b_3(t)=0$
for all $t\in I_{2n}=[t_{2n},t_{2n+1})$ and $b_1\in C^1(\bar I_{2n})$ for all $\ n\in\nat;$

\vskip+5 pt

 $jj_w$) $\vert b_3(t)\vert \le M_{2n+1}$,
  $b_1(t)=0$
for all $t\in I_{2n+1}=[t_{2n+1},t_{2n+2})$ and $b_3\in C(\bar I_{2n+1})$ for all  $\ n\in\nat.$

\vspace{0.5cm}
We emphasize that
 in the case without delay, since we deal only with the standard energy  $E_S(\cdot),$ the set $\tilde \omega$ where the negative damping is localized may be any subset of $\Omega,$ not necessarily a subset of $\omega.$

\begin{Remark}{\rm
Combining the results and the methods used so far, we can obtain stability results for problems with distributed or localized positive-negative damping with delay. We recall that the case of distributed positive-negative damping without delay was studied in \cite{Genni2}.
}
\end{Remark}
\subsection{Positive--negative damping without delay: the linear case}
\label{PNLinear}

\hspace{5mm}

In the linear case, as for the case with delay feedback, we can use a more explicit
observability constant in the interval $I_{2n}$ where only the positive damping is present.

Consider problem (\ref{cons1}) -- (\ref{cons2}) and assume that the observability inequality (\ref{OC}) holds (with now $W_1$ instead of $W$).

We can restate Proposition \ref{obsexplicit}.

\begin{Proposition}\label{obsexplicit2}
Assume $\mbox{\rm j)}$ and $f\equiv 0\,.$
Moreover,
we assume that the observability inequality $(\ref{OC})$ holds for every  time $T>\overline T$
and that, denoting  $T^*:=\inf_n \{T_{2n}\},$
\begin{equation}\label{ampiezze2}
T^*>\overline T.
\end{equation}
Then, for any solution of system $(\ref{4.1})-(\ref{4.2})$, we have
\begin{equation}\label{stimabuona2Anew2}
E_S(t_{2n+1})\le \hat c_n E_S(t_{2n}),\quad \forall\ n\in\nat,
\end{equation}
where
\begin{equation}\label{nnew2}
\displaystyle{
\hat c_n=
\frac{2c(1+4{C_1}^2T_{2n}^2M_{2n}^2) }
{m_{2n}+ 2c(1+4{C_1}^2T_{2n}^2M_{2n}^2)}
},
\end{equation}
  $c$ being the observability constant in $(\ref{OC})$ corresponding to the time $T^*$ and $C_1$ the constant in the norm embedding $(\ref{embeddi})$  between $W_1$ and $H$.
\end{Proposition}

Combining the previous proposition with (\ref{stimader2Cris3_senzaritardo1})
we can obtain the following stability result.

\begin{Theorem}\label{stab2Cris5_senzaritardo2}
Assume hypotheses of Proposition  $\ref{obsexplicit2},$
 $\mbox{\rm jj)}$, $(\ref{NL1})$ and $(\ref{length}).$
If
\begin{equation}\label{starstarNL_12}
\sum_{n=0}^\infty  [2C_3M_{2n+1}T_{2n+1}+\ln \hat c_n]=-\infty,
\end{equation}
where $\hat c_n$ is as in $(\ref{nnew2})$
and $C_3$ is the constant in the norm embedding $(\ref{embeddi})$ between $W_3$ and $H,$
then
system $(\ref{4.1})-(\ref{4.2})$ is asymptotically stable, that is any solution  $u$ of $(\ref{4.1})-(\ref{4.2})$ satisfies
$E_S(t)\rightarrow 0$ as $t\rightarrow +\infty.$
\end{Theorem}

\begin{Remark} {\rm
As in Remark \ref{remarklinear}, one can prove that
(\ref{starstarNL_12}) is satisfied if
$$
\sum_{n=0}^\infty M_{2n+1}T_{2n+1}<+\infty ,\quad \mbox{\rm and}\quad
\sum_{n=0}^\infty \frac {m_{2n}}
{1+4C_1^2T^2_{2n}M^2_{2n}}=+\infty.
$$}
\end{Remark}

\bigskip
\noindent{\bf Acknowledgements}
The authors are members of the Gruppo Nazionale per l'Analisi Matematica, la Probabilit\`a e le loro Applicazioni (GNAMPA) of the Istituto Nazionale di Alta Matematica (INdAM). Their
research is  partially supported by the GNAMPA
project  2015 {\it Analisi e controllo di equazioni a derivate parziali nonlineari}.


\begin{thebibliography}{10}
\bibitem{aloui}
F.~Aloui, I.~Ben Hassen and A.~Haraux.
\newblock Compactness of trajectories to some nonlinear second order evolution equations
and applications.
\newblock {\em J. Math. Pures Appl.}, 100: 295 -- 329, 2013.

\bibitem{aloui1}
F.~Aloui and A.~Haraux.
\newblock  Sharp ultimate bounds of solutions to a class of second
order linear evolution equations with bounded forcing.
\newblock {\em J. Funct. Anal.}, 265: 2204  -- 2225, 2013.

\bibitem{ANP2012}
K.~Ammari, S.~Nicaise and C.~Pignotti.
\newblock
Stabilization by switching time--delay.
\newblock {\em Asymptot. Anal.}, 83: 263--283, 2013.

\bibitem{BLR}
C.~Bardos, G.~Lebeau and J.~Rauch.
\newblock Sharp sufficient conditions for the observation, control and
stabilization of waves from the boundary.
\newblock {\em SIAM J. Control Optim.}, 30: 1024--1065, 1992.

\bibitem{lasiecka}
I.~Chueshov, M.~Eller and I.~ Lasiecka.
\newblock On the attractor for a semilinear wave equation with critical exponent and
nonlinear boundary dissipation.
\newblock {\em Comm. Partial Differential Equations}, 27: 1901–-1951, 2002.

\bibitem{Datko}
R.~Datko.
\newblock Not all feedback stabilized hyperbolic systems are robust
with respect to small time delays in their feedbacks.
\newblock {\em SIAM J. Control Optim.}, 26: 697--713, 1988.

\bibitem{DLP}
R.~Datko, J.~Lagnese and M.~P.~Polis.
\newblock An example on the effect of time delays in boundary feedback stabilization of
wave equations.
\newblock  {\em SIAM J. Control Optim.}, 24: 152--156, 1986.

\bibitem{ds}
\newblock  C.W.~De Silva.
 \newblock ``Vibration and Shock Handbook", \newblock Mechanical
Engineering, CRC Press 2005.

\bibitem{GenniSicon}
G.~Fragnelli and D.~Mugnai.
\newblock Stability of solutions for some classes of nonlinear damped wave equations
wave equations.
\newblock {\em SIAM J. Control Optim.}, 47: 2520--2539, 2008.

\bibitem{Genni2}
G.~Fragnelli and D.~Mugnai.
\newblock Stability of solutions for nonlinear wave equations with a positive--negative damping
\newblock {\em Discrete Contin. Dyn. Syst. Ser. S}, 4: 615--622, 2011.

\bibitem{GGH}
M.~Ghisi, M.~Gobbino and A.~Haraux.
\newblock The remarkable effectiveness of time-dependent damping terms for second order evolution equations.
\newblock {\em ArXiv:1506.06915}, 2015.


\bibitem{Gugat}
M.~Gugat.
\newblock  Boundary feedback stabilization by time delay for one-dimensional wave equations.
\newblock {\em IMA J. Math. Control Inform.}, 27: 189--203, 2010.


\bibitem{HMV}
A.~Haraux, P.~Martinez and J.~Vancostenoble.
\newblock     Asymptotic stability for intermittently controlled second--order
evolution equations.
\newblock {\em SIAM J. Control Optim.}, 43: 2089--2108, 2005.

\bibitem{KN}
\newblock  S.~Konabe and T.~Nikuni.
 \newblock Coarse--Grained Finite--Temperature
Theory for the Bose Condensate in Optical Lattices.
\newblock {\em J.
Low Temp. Phys.}, 150: 12--46, 2008.

\bibitem{LaMK}
\newblock  A.C.~Lazer and P.J.~McKenna.
\newblock Large--amplitude periodic
oscillations in suspension bridges: some new connections with
nonlinear analysis.
\newblock {\em SIAM Review}, 32: 537--578, 1990.


\bibitem{l}
H.A~Levine.
\newblock Some additional remarks on the nonexistence of global solutions to nonlinear wave equations.
\newblock {\em SIAM J. Math. Anal.}, 5: 138--146, 1974.

\bibitem{lps}
H.A~Levine, S.R.~Park and J.~Serrin.
\newblock Global existence and global nonexistence of solutions of the Cauchy problem for a nonlinearly damped wave equation.
\newblock {\em J. Math. Anal. Appl,}, 228: 181--205, 1998.


\bibitem{mm1}
\newblock A.~Marino and D.~Mugnai.
\newblock Asymptotically critical points and their
multiplicity.
\newblock{\em Topol. Methods Nonlinear Anal.} 19: 29--38, 2002.

\bibitem{mm2}
\newblock A.~Marino and D.~Mugnai.
\newblock Asymptotical multiplicity and some
reversed variational inequalities.
\newblock {\em Topol. Methods Nonlinear Anal.}20:  43--62, 2002.

\bibitem{bila}
D.~Mugnai.
\newblock On a "reversed" variational inequality,
\newblock {\em Topol. Methods Nonlinear Anal.}, 17: 321--358, 2001.


\bibitem{NPSicon06}
S.~Nicaise and C.~Pignotti.
\newblock Stability and instability results of the wave equation with a delay term in the boundary
or internal feedbacks.
\newblock {\em SIAM J. Control Optim.}, 45: 1561--1585, 2006.

\bibitem{ADE2012}
S.~Nicaise and C.~Pignotti.
\newblock Asymptotic stability of second--order evolution equations with
intermittent delay.
\newblock {\em Adv. Differential Equations}, 17: 879--902, 2012.


\bibitem{NP2014}
S.~Nicaise and C.~Pignotti.
\newblock
Stability results for second--order evolution equations with switching time--delay.
\newblock {\em J. Dynam. Differential Equations}, 26: 781--803, 2014.

\bibitem{Pazy}
A.~Pazy.
\newblock {\em Semigroups of linear operators and applications to partial differential equations}, Vol. 44 of {\em Applied Math. Sciences.} Springer-Verlag, New York, 1983.


\bibitem{Pignotti15}
C.~Pignotti.
\newblock
Stability results for second--order evolution equations with memory
and
switching time--delay.
\newblock {\em ArXiv:1507.03391}, 2015.

\bibitem{S}
\newblock  G.~Somieski.
\newblock Shimmy analysis of a simple aircraft nose
landing gear model using different mathematical methods.
 \newblock{\em
Aerosp. Sci. Technol.}1: 545--555, 1997



\bibitem{XYL}
G.~Q.~Xu, S.~P.~Yung and L.~K.~Li.
\newblock Stabilization of wave systems with input delay in the boundary control.
\newblock {\em ESAIM: Control Optim. Calc. Var.},
12: 770--785, 2006.



\bibitem{Zuazua}
E.~Zuazua.
\newblock Exponential decay for the semi-linear wave equation with locally distributed damping.
\newblock {\em Comm. Partial Differential Equations}, 15:205--235, 1990.
\end{thebibliography}
\end{document}